\documentclass[a4paper,11pt]{article}

\usepackage[utf8]{inputenc}
\usepackage[T1]{fontenc}
\usepackage{graphicx}
\usepackage{cancel}
\usepackage{amsmath,amsfonts,amssymb,amsthm}
\usepackage{mathtools}
\usepackage{mathrsfs}
\usepackage{xcolor}
\usepackage{hyperref}
\usepackage{booktabs}
\usepackage{bbold}
\usepackage{bm}
\usepackage{color}
\usepackage[font=footnotesize]{caption}
\usepackage{enumitem}
\usepackage{empheq}
\usepackage{framed}
\usepackage{fullpage}
\usepackage{hyperref}
\usepackage{soul}
\usepackage{dirtytalk}
\usepackage{pstricks}
\usepackage{appendix}
\usepackage{dsfont}
\usepackage{tcolorbox}
\usepackage{bbm}

\newcommand{\diag}{\operatorname{Diag}}
\newcommand{\off}{\operatorname{Off}}

\mathtoolsset{showonlyrefs}

\definecolor{myred}{HTML}{880000}
\definecolor{mygreen}{HTML}{008800}
\definecolor{myblue}{HTML}{000088}
\definecolor{linkblue}{HTML}{0000BB}

\hypersetup{
  colorlinks,
  linkcolor={myred},
  citecolor={myblue},
  urlcolor={linkblue}
}

\newcommand{\E}{{\mathbb E}}

\newcommand{\Var}{\operatorname{Var}}

\renewcommand{\leq}{\leqslant}
\renewcommand{\epsilon}{\varepsilon}

\renewcommand{\le}{\leqslant}
\renewcommand{\ge}{\geqslant}

\newcommand{\argmin}{\mathop{\mathrm{argmin}}}

\DeclareMathOperator{\tr}{Tr}

\usepackage{xspace}

\newtheorem{proposition}{Proposition}
\newtheorem{theorem}{Theorem}
\newtheorem*{theorem*}{Theorem}
\newtheorem{lemma}{Lemma}

\theoremstyle{definition}

\theoremstyle{remark}
\newtheorem{remark}{Remark}

\title{Covariance Estimation under Missing Observations and $L_4-L_2$ Moment Equivalence}
\author{
 Pedro Abdalla\thanks{Department of Mathematics, ETH Z\"urich,  \href{mailto:pedro.abdallateixeira@ifor.math.ethz.ch}{pedro.abdallateixeira@ifor.math.ethz.ch}}
}

\begin{document}

\maketitle

\begin{abstract}
We consider the problem of estimating the covariance matrix of a random vector by observing i.i.d samples, where entry of the sampled vector is missing with probability $p$. Under the standard $L_4-L_2$ moment equivalence assumption, we construct the first estimator that simultaneously achieves optimality with respect to the parameter $p$ and recovers the optimal convergence rate for the classical covariance estimation problem when $p=1$.
\end{abstract}

\section{Introduction}

High-dimensional covariance estimation is one of the most fundamental problems in the intersection of probability and statistics. On the applied side, it is a fundamental task for PCA or linear regression \cite{wainwright2019high}. On the theoretical side, the non-asymptotic properties of isotropic sample covariance matrices have been extensively studied \cite{adamczak2010quantitative,mendelson2012generic,vershynin2012close,tikhomirov2018sample,srivastava2013covariance,guedon2017interval} due to a famous question by Kannan, Lovász and Simonovits \cite{kannan1997random} and further generalized to the anisotropic case \cite{liaw2017simple,koltchinskii2017operators,abdalla2022}. Although the sample covariance matrix seems to be the most natural choice of estimator, its performance is suboptimal when the input data lacks a strong decay in the tail. Specifically, the convergence rate with respect to the confidence level $\delta$ is quite slow.

Motivated by this fact, a line of work in robust statistics, pioneered by Catoni \cite{catoni2012challenging}, studied the so-called sub-Gaussian estimators. These estimators are defined to be estimators that perform as good as the empirical mean under the Gaussian distribution. Many estimators have been proposed for the covariance estimation problem (see \cite{ke2019user} for a survey), in particular, there are now sub-Gaussian estimators under minimal assumptions on the data distribution \cite{abdalla2022,oliveira2022improved}.

On the other hand, data may be corrupted by noise. In \cite{lounici2014high}, Lounici addressed the so-called covariance estimation problem with \emph{missing observations}, motivated by applications in climate change, gene expression and cosmology. His work considers i.i.d observations, where each entry is ``missing" with probability $p$. We highlight that the missing observations model is a standard notion in the literature, extending beyond the covariance estimation setting, see \cite{elsener2019sparse,loh2011high} and the references therein.

The goal of this work is to design an estimator that simultaneously achieves the following properties:
\begin{itemize}
\item \textbf{Missing Observations:} We allow the data to have missing observations and heavy tails. We construct an estimator with minimax optimal convergence rate without assuming any knowledge of $p$. Remarkably, we show that dependence on $p$ is universal, meaning that it does not depend on the distribution of the data. 

\item \textbf{Dimension-Free:} The convergence rate scales with the effective rank $\mathbf{r}(\Sigma)$ rather than the dimension $d$,
\begin{equation*}
    \mathbf{r}(\Sigma):=\frac{\tr(\Sigma)}{\|\Sigma\|}.
\end{equation*}
This is an important aspect in high dimensional settings when the dimension $d$ is at least the sample size $N$. 

\item \textbf{Heavy-Tails:} We allow the distribution to have heavy tails, requiring only the existence of four moments satisfying minimal assumptions. Moreover, the result is as sharp as if the data were Gaussian (up to an absolute constant).

\end{itemize}

We begin with the rigorous definition of the model. We say that a centred random vector $X$ satisfies the $L_4-L_2$ \emph{moment equivalence} (\emph{hypercontractivity}) with constant $\kappa \ge 1$, if for all $v \in S^{d-1}$,
\begin{equation}
\label{eq:momeqv}
(\E \langle X, v\rangle^{4})^{1/4} \le \kappa(\E \langle X, v\rangle^{2})^{1/2}.
\end{equation}
Here we always assume that the data satisfies the $L_4-L_2$ moment equivalence with an absolute constant $\kappa>0$, i.e, the constant $\kappa$ is a fixed real number that does not depend on any other parameter. A vast class of distributions satisfies the moment equivalence assumption mentioned above, with $\kappa$ being a small absolute constant. Examples include sub-gaussian random vectors, sub-exponential random vectors with bounded $\psi_\alpha$ norm, as well $t$-student distributions with a sufficiently large degree of freedom \cite{mendelson2020robust}.
\

We say that the sample $Y_1,\ldots,Y_N$ is $p$-\emph{sparsified} if it is obtained from the sample $X_1, \ldots X_N$ of independent copies of $X$ by multiplying each entry of the $X_i$'s by an independent 0/1 Bernoulli random variable with mean $p$. In concise terminology, we say that the data is sampled from $X\odot \mathbf{p}$, where $\mathbf{p} \in \{0,1\}^d$ is a random vector with i.i.d entries Bernoulli $p$, and the notation $\odot$ simply denotes the standard entrywise product. The choice of zero to represent missing information is merely for convenience and could be replaced by any other value. Now, we present the main result of this manuscript.
\begin{theorem}[Main result]
\label{main_result}
Assume that $X$ is a zero mean random vector in $\mathbb{R}^d$ with covariance matrix $\Sigma$ satisfying the $L_4-L_2$ moment equivalence assumption with an absolute constant $\kappa$. Fix the confidence level $\delta \in (0,1)$. Suppose that $Y_1,\ldots,Y_N$ are i.i.d samples distributed as $X\odot p$, where $\mathbf{p}=(p_1,\ldots,p_d) \in \mathbb{R}^d$ is a random vector with  i.i.d Bernoulli entries with parameter $p$. Then there exists an estimator $\widehat{\Sigma}(N,\delta)$ depending only on the sample $Y_1,\ldots,Y_N$ and $\delta$ satisfying that, with probability at least $1-\delta$, 
\begin{equation*}
\|\widehat{\Sigma} - \Sigma \| \leq \frac{C(\kappa)}{p} \|\Sigma\| \left(\sqrt{\frac{r(\Sigma)+\log(1/\delta)}{N}}\right).
\end{equation*}
Here $C(\kappa)>0$ is an absolute constant depending only on $\kappa$.
\end{theorem}

\paragraph{Literature Review:}
We remark that several results for covariance estimation under missing observations were obtained in the literature, for example \cite{lounici2014high,Rao,pavez2020covariance,park2021,klochkov2020uniform,cai2016minimax}. However, none of the previous results has been able to simultaneously scale correctly with the factor of $p$ and recover a sub-Gaussian estimator when $p=1$, as established in \cite{abdalla2022,oliveira2022improved}, even when the data is Gaussian. Moreover, the convergence rate is optimal up to an absolute constant: When $p=1$, a classical result by Lounici and Koltchinskii \cite[Theorem 4]{koltchinskii2017operators} states that if $G_1,\ldots, G_N$ are i.i.d mean zero Gaussian vectors with covariance matrix $\Sigma$ and $N\ge r(\Sigma)$, then
\begin{equation*}
c\|\Sigma\| \sqrt{\frac{r(\Sigma)}{N}}\le \mathbb{E}\left \|\frac{1}{N}\sum_{i=1}^N G_i\otimes G_i - \Sigma \right\| \le  C \|\Sigma\| \sqrt{\frac{r(\Sigma)}{N}}.
\end{equation*}
This essentially shows optimality with respect to the effective rank, as we expect that the empirical covariance for Gaussian distributions to be sharp in expectation. They also showed that that the expectation is tightly concentrated around the mean, and our quantitative convergence rate with respect to $\delta$ matches their result up to an absolute constant. Both are indeed optimal among all (measurable) estimators for the covariance; see \cite{mendelson2020robust,abdalla2022} for a more technical discussion. Intuitively, it should be not surprising that we cannot beat the Gaussian decay.

In addition to this, the dependence with respect to $p$ is also optimal thanks to a minimax lower bound from Lounici \cite[Theorem 2]{lounici2014high}. In a nutshell, his result shows that there exists absolute constants $c_1,c_2>0$ for which
\begin{equation*}
    \inf_{\widehat{\Sigma}} \sup_{\mathbb{P}} \mathbb{P}\left(\|\widehat{\Sigma}-\Sigma\|\ge \frac{c_1}{p}\|\Sigma\|\sqrt{\frac{r(\Sigma)}{N}}\right) \ge c_2.
\end{equation*}
Here, the infimum is taken with respect to all estimators that depend only on the data, and the supremum is taken over all possible distributions with covariance matrix $\Sigma$. This implies that our main result captures the optimal dependence with respect to $p$ as well.

It is important to note that the main drawback of our result is that our estimator is not computationally tractable. The primary focus of this work is on the information-theoretic limits of covariance estimators, specifically our main contribution is demonstrating the possibility of constructing an optimal data-driven estimator for covariance under minimal assumptions on the data (albeit it is computationally infeasible).

To the best of our knowledge, there are no computable estimators for the covariance matrix under heavy tails, even in the case without missing observations. We leave it as an important open problem.

\paragraph{Proposed Estimator:}
The startpoint to construct our estimator is the following observation: The expectation of the covariance matrix of $Y$ scales differently for the diagonal part and the off-diagonal part of its covariance matrix. More accurately,
\begin{equation*}
    \mathbb{E}Y\otimes Y = p\diag(\Sigma) + p^2\off(\Sigma).
\end{equation*}
We can ``invert" the equality above to get the dependence between the true covariance and the data, namely
\begin{equation*}
\Sigma = p^{-1}\diag(\mathbb{E}Y\otimes Y) + p^{-2}\off(\mathbb{E}Y\otimes Y).
\end{equation*}
A natural approach would be to replace the unknown term $\mathbb{E}Y\otimes Y$ by its sample covariance, but this is not enough when we consider heavy tailed data $X_1,\ldots,X_N$, as discussed above. In fact, we define the truncation function
\begin{equation}
\label{eq:truncfunction}
\psi(x) = 
    \begin{cases}
      x,\quad &\textrm{for}\; x \in [-1, 1],
      \\
      \operatorname{sign}(x),\quad &\textrm{for}\; |x| > 1,
    \end{cases}
\end{equation}
to robustify our estimator in each direction of the sphere. The idea here is to estimate the matrix through its quadratic form. Next, we describe the estimator's final form. We estimate the diagonal and off-diagonal part separately,
\begin{equation*}
\widehat{\Sigma_1}(\lambda_1):=\argmin_{\Sigma_1\in \mathbb{S}_{+}^d|\off(\Sigma_1)=0} \sup_{v\in S^{d-1}} |v^T \Sigma_1 v - \frac{1}{n\lambda_1}\sum_{i=1}^n \psi(\lambda_1v^T\diag(Y_i\otimes Y_i)v)|,
\end{equation*}
\begin{equation*}
\widehat{\Sigma_2}(\lambda_2):=\argmin_{\Sigma_2\in \mathbb{S}_{+}^d|\diag(\Sigma_2)=0} \sup_{v\in S^{d-1}} |v^T \Sigma_2 v - \frac{1}{n\lambda_2}\sum_{i=1}^n \psi(\lambda_2v^T\off(Y_i\otimes Y_i)v)|,
\end{equation*}
where $\mathbb{S}_{+}^d$ is the set of $d$ by $d$ positive semi-definite matrices and $S^{d - 1}$ denotes the unit sphere in $\mathbb{R}^d$, and the final estimator becomes 
\begin{equation*}
\widehat{\Sigma} = \frac{1}{\widehat{p}}\diag(\widehat{\Sigma_1}) + \frac{1}{\widehat{p}^2}\off(\widehat{\Sigma_2}).
\end{equation*}
Here, $\widehat{p}$ is an estimator for the parameter $p$, and the choice of the truncation levels $\lambda_1,\lambda_2$ will be clarified in what follows.

As mentioned before, the main drawback of our estimation is that it is not computational tractable. Indeed, $\widehat{\Sigma}_1,\widehat{\Sigma}_2$ does not seem to be computable in polynomial time as a (sub)-gradient descent/ascent type method might get stuck in local optima, and analyzing it is out of the scope of this text. We remark that a similar optimization problem appears in \cite{depersin2022optimal}, with the quadratic forms replaced by linear forms. Unfortunately, even in that case it is an open problem to come up with an (analyzable) algorithm.

From a more practical perspective, it might be possible that, under some stronger concentration assumption on the data, we can avoid evaluating the truncation function $\psi$ in each direction $v$ on the Euclidean sphere, and replace the supremum over the Euclidean sphere in the definition of $\widehat{\Sigma}_1,\widehat{\Sigma}_2$ by other (more tractable) quantity. This would make the optimization problem much easier to be solved in polynomial time. We leave this as an interesting question to pursue in a future work.

The construction of the estimator and its analysis share similarities with the "trimmead covariance" estimator proposed by Zhivotovskiy and the author \cite{abdalla2022}. However, we need to break into diagonal and off-diagonal parts to take in account the different scales with $p$. Indeed, the main technical difficulty arises in controlling the random quadratic form to get the optimal dependence with respect to $p$, mainly in the off-diagonal case. A direct approach faces the difficulty that we no longer have a positive semidefinite matrix, making it challenging to capture cancellations. Conversely, an indirect approach, expressing the off-diagonal part as the total part minus the diagonal part, leads to sub-optimality with respect to $p$. Thus, we need to carefully balance these two approaches.

\paragraph{Organization.}
The rest of the paper is organized as follows: In Section \ref{sec:oracle}, we assume the knowledge of certain parameters to simplify the analysis of the estimator and derive sharp convergence rates. We then systematically relax these assumptions in Section \ref{sec:proof_main} by estimating each parameter separately in individual subsections. The last subsection of Section \ref{sec:proof_main} is devoted to the formal construction of the estimator and the proof of the main result.
\paragraph{Notation.}
Throughout this text $C,c>0$ denote an absolute constant that may change from line to line. For an integer $N$, we set $[N] = \{1, \ldots, N\}$. For any two functions (or random variables) $f, g$ defined in some common domain, the notation $f \lesssim g$ means that there is an absolute constant $c$ such that $f \le cg$ and $f\sim g$ means that $f\lesssim g$ and $g\lesssim f$. Let $\mathbb{S}_{+}^d$ denote the set of $d$ by $d$ positive-definite matrices. The symbols $\|\cdot\|,\|\cdot\|_{F}$ denote the operator norm and the Frobenius norm of a matrix, respectively. Let $ \mathcal{KL}(\rho, \mu) = \int\log\left(\frac{d\rho}{d\mu}\right)d\rho$
denote the Kullback-Leibler divergence between a pair of measures $\rho$ and $\mu$. We write $\rho \ll \mu$ to indicate that the measure $\rho$ is absolutely continuous with respect to the measure $\mu$. For a vector $X\in \mathbb{R}^d$, the tensor product $\otimes$ is defined as $X\otimes X:= XX^T \in \mathbb{R}^{d\times d}$.

\section{Oracle Estimator}
\label{sec:oracle}
In this section, we prove our main result under the assumption that we know the effective rank of the covariance matrix $r(\Sigma)$, the trace of the covariance matrix $\tr(\Sigma)$, and the sparsifying factor $p$. These assumptions will be further relaxed in the next section. Our main goal is to prove the following result.
\begin{proposition}
\label{oracle_estimator}
Assume that $X$ is a mean zero random vector in $\mathbb{R}^d$ with covariance matrix $\Sigma$ satisfying the $L_4-L_2$ moment equivalence assumption. Fix the confidence level $\delta \in (0,1)$. Suppose that $Y_1,\ldots,Y_N$ are i.i.d samples from $X\odot \mathbf{p}$. Then there exists $\lambda_1,\lambda_2>0$ depending only on $\tr(\Sigma),\|\Sigma\|$ and $p$ for which, with probability at least $1-\delta$,
\begin{equation*}
\begin{split}
  &\max\{\|p^{-1}\widehat{\Sigma_1}(\lambda_1) - \diag(\Sigma)\|, \|p^{-2}\widehat{\Sigma_2}(\lambda_2)-\off(\Sigma)\| \} \\
  &\le \frac{C(\kappa)}{p}\|\Sigma\|\left(\sqrt{\frac{r(\Sigma)+\log(1/\delta)}{N}}\right).
  \end{split}
\end{equation*}
Here $C(\kappa)>0$ is an absolute constant depending only on $\kappa$.
\end{proposition}

Our analysis is based on the variational principle pioneered by O. Catoni \cite{catoni2012challenging,catoni2007pacbayes,catoni2018dimension} and further developed in many applications related to high dimensional probability and statistics \cite{catoni2007pacbayes,catoni2016pac,catoni2017dimension,catoni2018dimension,zhivotovskiy2021dimension,mourtada2019exact}. In most of the applications of the variational principle, the following lemma serves as a key stepping stone.
\begin{lemma}
\label{lem:pacbayes}
Assume that $X_i$ are i.i.d. random variables defined on some measurable space. Let $\Theta$ be a subset of $\mathbb{R}^p$ for some $p \ge 1$, $\mu$ be a a fixed distribution on $\Theta$, and $\rho$ be any distribution on $\Theta$ satisfying that $\rho \ll \mu$. Then, simultaneously for any such $\rho$, with probability at least $1 - \delta$,
\[
\frac{1}{N}\sum\limits_{i = 1}^N\mathbb{E}_{\rho}f(X_i, \theta) \le \mathbb{E}_{\rho}\log(\mathbb{E}_X e^{f(X, \theta)}) + \frac{\mathcal{KL}(\rho, \mu) + \log(1/\delta)}{N}.
\]
Here $\theta$ is distributed according to $\rho$. 
\end{lemma}
The proof can be found in \cite{catoni2007pacbayes,zhivotovskiy2021dimension} and will be omitted. The next lemma is a technical fact that allow to "convexify" the truncation function $\psi$. Indeed, it is easy to see that the function $e^\psi(x)$ is bounded by $(1+x+x^2)$ that still not convex, but if we add a suitable quadratic term, then it becomes convex.
\begin{lemma}
\label{lem:almostconvex}
Let $\psi$ be the truncation function from \eqref{eq:truncfunction}, and let $Z$ be a random variable with finite second moment. Then the following holds
\[
\psi(\mathbb{E} Z) \le \mathbb{E}\log(1 + Z + Z^2) + \min\left\{1,  \frac{\mathbb{E} Z^2}{6}\right\}.
\]
Moreover, for any $a > 0$, 
\begin{equation*}
\begin{split}
 \mathbb{E}\log(1 + Z + Z^2) &+ a \min\left\{1,  \frac{\mathbb{E} Z^2}{6}\right\}\\
 &\le \mathbb{E}\log\left(1 + Z + \left(1 + \frac{(7 + \sqrt{6})(\exp(a) - 1)}{6}\right) Z^2\right).
\end{split}
\end{equation*}
\end{lemma}
This result has been previously used in \cite{catoni2007pacbayes,zhivotovskiy2021dimension,abdalla2022}. For the sake of completeness, we include a proof at the end of this section. Now we start with the facts specifically derived for the missing observation case. The next result is crucial to establish the right dependence on $p$. The proof is deferred to the end of this section.
\begin{lemma}
\label{lemma:diagxoff}
Let $Y$ as above. For every $v\in S^{d-1}$, we have
\begin{equation*}
    \mathbb{E} (v^T\diag(Y\otimes Y)v)^2 \le 2p\kappa^4\|\diag(\Sigma)\|^2
\end{equation*}
and 
\begin{equation*}
    \mathbb{E} (v^T\off(Y\otimes Y)v)^2 \le 4p^2\kappa^4\|\Sigma\|^2.
\end{equation*}
\end{lemma}
The main idea behind the proof of Proposition \ref{oracle_estimator} consists in using the variational principle twice, one for the diagonal part and the other for the more delicate off-diagonal part. 
\begin{proof}
\textbf{Diagonal Part}:
We start by defining the parameter space of interest, namely
\[
\Theta = \mathbb{R}^d \times \mathbb{R}^d.
\]
Choose $\mu$ to be a product of two zero mean multivariate Gaussians with mean zero and covariance $\beta^{-1}I_d$, where $\beta>0$ will be chosen later. For each $v \in S^{d - 1}$, let $\rho_{v}$ be the product of two multivariate Gaussian distribution with mean $v$ and covariance $\beta^{-1}I_d$. By construction, $(\theta, \nu)$ is distributed according to $\rho_{v}$, therefore it satisfies that $\mathbb{E}_{\rho_{v}}(\theta, \nu) = (v,v)$. The standard formula for the $\mathcal{KL}$-divergence between two Gaussian measures \cite{pardo2018} implies that
\[
\mathcal{KL}(\rho_{v}, \mu) = \beta.
\]
Let $\lambda_1 > 0$ be a free parameter to be optimized later. By the first part of Lemma \ref{lem:almostconvex}, we have 
\begin{equation*}
\begin{split}
&\psi\left(\lambda_1 v^T\diag(Y\otimes Y)v\right) = \psi\left(\lambda_1 \mathbb{E}_{\rho_v}\theta^T\diag(Y\otimes Y)\nu\right)\\
&\le \mathbb{E}_{\rho_v}\log(1+\lambda_1 \theta^T\diag(Y\otimes Y)\nu+\lambda_1^2 (\theta^T\diag(Y\otimes Y)\nu)^2 ) + R.
\end{split}
\end{equation*}
where $R:=\min\{1,\lambda_1^2\mathbb{E}_{\rho_v}(\theta^T\diag(Y\otimes Y)\nu)^2/6\}$. Notice that $\mathbb{E}\theta_i^2 = \beta^{-1}+v_i^2$ and $\mathbb{E}\theta_i\theta_j = v_iv_j$ for all $i\in [d]$, therefore
\begin{equation*}
\begin{split}
&\mathbb{E}_{\rho_v}(\theta^T\diag(Y\otimes Y)\nu)^2 = \mathbb{E}_{\rho_v}\left(\sum_{i=1}^d\langle Y,e_i\rangle^2\theta_i\nu_i\right)^2 \\
&= \beta^{-2}\sum_{i=1}^d \langle Y,e_i\rangle^4 + \sum_{i,j=1}^d\langle Y,e_i\rangle^2\langle Y,e_j\rangle^2 v_i^2v_j^2 + 2\beta^{-1}\sum_{i=1}^d\langle Y,e_i\rangle^4 v_i^2 \\
&= \beta^{-2}\|\diag(Y\otimes Y)\|_F^2 + (v^T\diag(Y\otimes Y)v)^2 + 2\beta^{-1} \|\diag(Y\otimes Y)v\|_2^2,
\end{split}
\end{equation*}
By symmetry, $\mathbb{P}\left(\theta^T\diag(Y\otimes Y)\nu\ge v^T\diag(Y\otimes Y)v\right) \ge \frac{1}{4}$. To see this, observe that it is equal to
\begin{equation*}
\mathbb{P}(\langle \diag(Y\otimes Y)\theta,(\nu-v)\rangle +\langle \diag(Y\otimes Y)v,(\theta-v)\rangle \ge 0).
\end{equation*}
The second term is positive with probability one half. Conditioned on this event, the first term is positive with probability one half, and it is independence from the first term. We obtain that the probability of both are positive is at least one quarter. Therefore,
\begin{equation*}
\min\left\{1,\frac{\lambda_1^2}{6}(v^T\diag(Y\otimes Y) v)^2\right\} \leq 4 \mathbb{E}_{\rho_v}\min\left\{1,\frac{\lambda^2}{6}(\theta^T\diag(Y\otimes Y)\nu)^2\right\} .
\end{equation*}
By the second part of Lemma \ref{lem:almostconvex}, we have
\begin{equation*}
\begin{split}
&\psi(\lambda_1v^T\diag(Y\otimes Y)v)\\
&\leq \mathbb{E}_{\rho_v} \log (1+ \lambda \theta^T \diag(Y\otimes Y)\nu + C_1\lambda^2(\theta^T\diag(Y\otimes Y)\nu)^2) + R(Y,\beta),
\end{split}
\end{equation*}
where $R(Y,\beta):= \min\{1,2\lambda_1^2\beta^{-1}\|\diag(Y\otimes Y)v\|_2^2/6\}+\min\{1,\lambda_1^2\beta^{-2}\|Y\otimes Y\|_F^2/6\}$. For instance, let us focus on the first term. The goal is to apply Lemma \ref{lem:pacbayes} to the function $f$ defined below
\begin{equation*}
f(Y,\theta,\nu):= \log(1+\lambda_1 \theta^T\diag(Y\otimes Y)\nu + C_1\lambda_1^2(\theta^{T}\diag(Y\otimes Y)\nu)^2).
\end{equation*}

\noindent Using the numeric inequality $\log(1 + y) \le y$, valid for all $y \ge - 1$, followed by Fubini's theorem and Lemma \ref{lemma:diagxoff}, we have
\begin{align*}
&\mathbb{E}_{\rho_{v}}\log\mathbb{E} \left(1+\lambda_1 \theta^T\diag(Y\otimes Y)\nu + C_1\lambda_1^2(\theta^{T}\diag(Y\otimes Y)\nu)^2\right)
\\
&\le  \mathbb{E}_{\rho_{v}}\mathbb{E} \left(\lambda_1 \theta^T\diag(Y\otimes Y)\nu + C_1\lambda_1^2(\theta^{T}\diag(Y\otimes Y)\nu)^2\right)
\\
&\le p \lambda_1 v^T\diag(\Sigma)v + C_1\lambda_1^2 (p\beta^{-2}\kappa^4\tr^2(\Sigma) + 2\beta^{-1}p\kappa^4 \|\diag\Sigma\|^2 + p\kappa^4\|\Sigma\|^2 )\\
&\le p \lambda_1 v^T\diag(\Sigma)v + C_1\lambda_1^2 (p\beta^{-2}\kappa^4\tr^2(\Sigma) + 2\beta^{-1}p\kappa^4 \|\Sigma\|^2 + p\kappa^4\|\Sigma\|^2 ).
\end{align*}
Next, setting $\beta:=r(\Sigma)$ (which is at least one) and applying Lemma \ref{lem:pacbayes}, it follows that with probability at least $1-\delta$, for all $v\in S^{d-1}$,
\begin{equation*}
\begin{split}
&\frac{1}{N\lambda_1}\sum_{i=1}^N \psi(\lambda_1 v^T\diag(Y\otimes Y)v) \\
&\le p v^T\diag(\Sigma)v + C\lambda_1 p\|\Sigma\|^2\kappa^4 + \sum_{i=1}^n\frac{R_i}{N\lambda_1} + \frac{r(\Sigma)+\log(1/\delta)}{\lambda_1N}.
\end{split}
\end{equation*}
Here $R_i$ is an independent copy of $R$. We proceed to estimate the third term in the right-hand side. Clearly, since $\min\{1,2\lambda_1^2\beta^{-1}\|\diag(Y\otimes Y)v\|_2^2/6\}$ is bounded by one, its variance is bounded by its expectation. Therefore, by Bernstein's inequality it follows that with probability $1-\delta$,
\begin{equation*}
\begin{split}
&\frac{1}{\lambda_1 N}\sum_{i=1}^N\min\{1, \beta^{-2}\lambda_1^2\|Y\otimes Y\|_F^2/6\} \lesssim \mathbb{E} \beta^{-2}\|Y\otimes Y\|_F^2 + \frac{\log(1/\delta)}{\lambda_1 N}\\
&\lesssim \lambda_1p\kappa^4\|\Sigma\|^2 + \frac{\log(1/\delta)}{\lambda_1 N}.
\end{split}
\end{equation*}
An analogous computation shows the same estimate holds (up to an absolute constant) for the term $\min\{1,2\beta^{-1}\|\diag(Y\otimes Y)v\|_2^2/6\}$. Finally we conclude that, there exists an absolute constant $C>0$ such that, with probability at least $1-\delta$,
\begin{equation*}
\begin{split}
&\frac{1}{N\lambda_1}\sum_{i=1}^N \psi(\lambda_1 v^T\diag(Y_i\otimes Y_i)v)\\
&\le p v^T\diag(\Sigma)v + C\left(\lambda_1 p\|\Sigma\|^2\kappa^4 + \frac{r(\Sigma)+\log(1/\delta)}{\lambda_1N}\right).
\end{split}
\end{equation*}
We optimize the right-hand side over $\lambda_1>0$. More accurately, setting
\begin{equation*}
\lambda_1 =  \frac{1}{\|\Sigma\|\kappa^2p}\sqrt{\frac{r(\Sigma)+\log(1/\delta)}{N}}
\end{equation*}
we obtain that, with probability at least $1-\delta$,
\begin{equation*}
\frac{1}{N\lambda_1}\sum_{i=1}^N \psi(\lambda_1 v^T\diag(Y_i\otimes Y_i)v) \le p v^T\diag(\Sigma)v + C\kappa^2\sqrt{p}\sqrt{\frac{r(\Sigma)+\log(1/\delta)}{N}}.
\end{equation*}
We repeat the same arguments above for $\rho_{2,v}$ being a product measure between two Gaussians $\theta \sim N(v,\beta^{-1}I_d)$ and $\nu \sim N(-v,\beta^{-1}I_d)$. The argument follows exactly the same steps because $\psi$ is symmetric. Therefore, it also holds that with probability $1-\delta$,
\begin{equation*}
-\frac{1}{N\lambda_1}\sum_{i=1}^N \psi(\lambda_1v^T\diag(Y_i\otimes Y_i)v) \le -p v^T\diag(\Sigma)v + C\kappa^2\sqrt{p}\sqrt{\frac{r(\Sigma)+\log(1/\delta)}{N}}.
\end{equation*}
By union bound, we obtain a two-sided bound: With probability at least $1-\delta$,
\begin{equation*}
\|p^{-1}\widehat{\Sigma}_1 - \diag(\Sigma)\| \lesssim \frac{1}{\sqrt{p}}\kappa^2\|\Sigma\| \sqrt{\frac{r(\Sigma)+\log(1/\delta)}{N}}.
\end{equation*}
\textbf{Off-diagonal part}: We now proceed to the second part of the proof to deal with the off-diagonal part. We choose $\mu$ and $\rho(v)$ as before, and write
\begin{equation*}
\begin{split}
&\psi\left(\lambda_2 v^T\off(Y\otimes Y)v\right) = \psi\left(\lambda_2 \mathbb{E}_{\rho_v}\theta^T\off(Y\otimes Y)\nu\right)\\
&\le \mathbb{E}_{\rho_v}\log(1+\lambda_2 \theta^T\off(Y\otimes Y)\nu+\lambda_2^2 (\theta^T\off(Y\otimes Y)\nu)^2 ) + R.
\end{split}
\end{equation*}
where $R_2:=\min\{1,\lambda_2^2 \mathbb{E}_{\rho_v}(\theta^T \off(Y\otimes Y)\nu)^2 /6\}$. We have to deal with the quadratic form of the off-diagonal that requires a more delicate analysis. In fact,
\begin{equation*}
\mathbb{E}_{\rho_v}(\theta^T \off(Y\otimes Y)\nu)^2 = \mathbb{E}_{\rho_v}\sum_{i\neq j;k\neq l} \langle Y,e_i\rangle\langle Y,e_j\rangle \langle Y,e_k\rangle \langle Y,e_l\rangle \theta_i \nu_j \theta_k \nu_l
\end{equation*}
By independence between $\theta$ and $\nu$, it remains to analyze the term $\mathbb{E}\theta_i\theta_k\mathbb{E}\nu_j\nu_l$. We split the analysis in three cases: The first one when $k=i$ and $j=l$, the second when either $k=i$ and $j\neq l$ or $k\neq i$ and $j=l$, and finally the third one when $k\neq i$ and $j\neq l$. In the first case, the summation becomes 
\begin{equation*}
\sum_{i\neq j}\langle Y,e_i\rangle^2\langle Y,e_j\rangle^2(\beta^{-1}+v_i^2)(\beta^{-1}+v_j^2).
\end{equation*}
In the second case, the summation becomes
\begin{equation*}
\sum_{i\neq j\neq l} \langle Y,e_i\rangle^2 \langle Y,e_j\rangle \langle Y,e_l\rangle (\beta^{-1}+v_i^2)v_jv_l + \sum_{i\neq j\neq k} \langle Y,e_i\rangle \langle Y,e_j\rangle^2 \langle Y,e_k\rangle (\beta^{-1}+v_j^2)v_iv_k.
\end{equation*}
The third case is simpler,
\begin{equation*}
\sum_{i\neq j\neq k\neq l}  \langle Y,e_i\rangle \langle Y,e_j\rangle \langle Y,e_k\rangle \langle Y,e_l\rangle v_iv_jv_kv_l.
\end{equation*}
Observe that summing all terms that do not contain any $\beta$ factor, we obtain $(v^T\off(Y\otimes Y) v)^2$. As before, the goal is to apply Lemma \ref{lem:pacbayes} to the function $f$,
\begin{equation*}
f(Y,\theta,\nu):= \log(1+\lambda_2 \theta^T\off(Y\otimes Y)\nu + C_2\lambda_2^2(\theta^{T}\off(Y\otimes Y)\nu)^2),
\end{equation*}
where $C_2>0$ is a sufficiently large absolute constant.
\noindent Using again the numeric inequality $\log(1 + y) \le y$, Fubini's theorem, and Lemma \ref{lemma:diagxoff},  we have
\begin{align*}
&\mathbb{E}_{\rho_{v}}\log\mathbb{E} \left(1+\lambda_2 \theta^T\off(Y\otimes Y)\nu + C_2\lambda_2^2(\theta^{T}\off(Y\otimes Y)\nu)^2\right)
\\
&\le  \mathbb{E}_{\rho_{v}}\mathbb{E} \left(\lambda_2 \theta^T\off(Y\otimes Y)\nu + C_2\lambda_2^2(\theta^{T}\off(Y\otimes Y)\nu)^2\right).
\end{align*}
The first term is equal to $p^2\lambda_2v^T\off(\Sigma)v$. We know that all terms in the expansion of $\theta^T \off(Y\otimes Y) \nu$ without a $\beta$ factor add up $v^T \off(Y\otimes Y) v$ and its expectation is at most $4p^2\kappa^4\|\Sigma\|^2$ by Lemma \ref{lemma:diagxoff}. Next, we estimate the terms containing $\beta$ systematically. Using Cauchy-Schwarz inequality together with the moment equivalence for $X$, we obtain
\begin{equation*}
\begin{split}
&\beta^{-2}\sum_{i\neq j}\mathbb{E}\langle Y,e_i\rangle^2\langle Y,e_j\rangle^2  = \beta^{-2}p^2\sum_{i\neq j}\mathbb{E}\langle X,e_i\rangle^2\langle X,e_j\rangle^2 \le p^2\beta^{-2}\kappa^4\sum_{i\neq j}\Sigma_{ii}\Sigma_{jj}\\
& \le p^2\beta^{-2}\kappa^4 \tr^2(\Sigma).
\end{split}
\end{equation*}
Similarly, we obtain
\begin{equation*}
\mathbb{E}\sum_{i\neq j}\langle Y,e_i\rangle^2\langle Y,e_j\rangle^2(\beta^{-1}+v_i^2)(\beta^{-1}+v_j^2) \lesssim p^2\kappa^4\beta^{-2}\tr^2(\Sigma) + \beta^{-1}p^2\kappa^4\|\Sigma\|\tr(\Sigma).
\end{equation*}
It remains to analyze
\begin{equation*}
\mathbb{E} \sum_{i\neq j\neq l} \langle Y,e_i\rangle^2 \langle Y,e_j\rangle \langle Y,e_l\rangle \beta^{-1}v_jv_l + \sum_{i\neq j\neq k} \langle Y,e_i\rangle \langle Y,e_j\rangle^2 \langle Y,e_k\rangle \beta^{-1}v_iv_k.
\end{equation*}
We estimate the first term on the right-hand side as the second term is identically distributed. To this end, we apply H\"older's inequality with conjugate exponents $4/3$ and $4$, and the moment equivalence to obtain that
\begin{equation*}
\begin{split}
&\mathbb{E} \sum_{i\neq j\neq l} \langle Y,e_i\rangle^2 \langle Y,e_j\rangle \langle Y,e_l\rangle \beta^{-1}v_jv_l \le p^3 \sum_{i\neq j\neq l}\mathbb{E} \langle X,e_i\rangle^2 \langle X,e_j\rangle \langle X,e_l\rangle \beta^{-1}v_jv_l \\
&\le p^3\sum_{i;j\neq l}\mathbb{E}\langle X,e_i\rangle^2\langle X,e_j\rangle \langle X,e_l\rangle\beta^{-1}v_jv_l + \left|2p^3\sum_{j\neq l}\mathbb{E}\langle X,e_j\rangle^3\langle X,e_l\rangle\beta^{-1}v_jv_l \right|\\
&\lesssim \frac{p^3}{\beta} \left(\mathbb{E} \left[ (v^T\off(X\otimes X)v)\left(\sum_{i=1}^d \langle X,e_i\rangle^2\right)\right] + \kappa^4\sum_{j\neq l}(\Sigma_{ll})^{1/2} (\Sigma_{jj})^{3/2}v_lv_j\right)\\
&\lesssim p^3\left[ (v^T\off(X\otimes X)v)\left(\sum_{i=1}^d \beta^{-1}\langle X,e_i\rangle^2\right)\right] + p^3\beta^{-1}\kappa^4\|\Sigma\|\tr(\Sigma)\\
& \le p^3 \mathbb{E}(v^TX\otimes X v - v^T\diag(X\otimes X)v)(\beta^{-1}\tr(X\otimes X)) + 2p^3\frac{\kappa^4}{\beta}\|\Sigma\|\tr(\Sigma)\\
& \le p^3\mathbb{E} \langle X,v\rangle^2\beta^{-1}\tr(X\otimes X) + 2p^3\beta^{-1}\kappa^4\|\Sigma\|\tr(\Sigma)\\
&\lesssim p^3 \kappa^4( \|\Sigma\|^2 + \beta^{-2}\tr^2(\Sigma) + \beta^{-1}\|\Sigma\|\tr(\Sigma)),
\end{split}
\end{equation*}
where the last inequality follows from the arithmetic-geometric inequality. Putting all together, we conclude that 
\begin{equation*}
\mathbb{E}\mathbb{E}_{\rho_v}(\theta^T\off(Y\otimes Y) \nu)^2 \lesssim p^2 \beta^{-2}\kappa^4 \tr^2(\Sigma) + p^2\kappa^4\|\Sigma\|^2 + p^3\beta^{-1}\kappa^4 \tr(\Sigma)\|\Sigma\|.
\end{equation*}
Therefore, setting $\beta = r(\Sigma)$, it follows that
\begin{align*}
&\mathbb{E}_{\rho_{v}}\log\mathbb{E} \left(1+\lambda_2 \theta^T\diag(Y\otimes Y)\nu + C_2\lambda_2^2(\theta^{T}\diag(Y\otimes Y)\nu)^2\right)
\\
&\le \lambda_2v^T\off(Y\otimes Y) v + C_2 \lambda_2^2\delta^2\kappa^4\|\Sigma\|^2.
\end{align*}
Finally we conclude that there exists an absolute constant $C_2'>0$ for which, with probability at least $1-\delta$,
\begin{equation*}
\begin{split}
&\frac{1}{N\lambda_2}\sum_{i=1}^N \psi(\lambda_2 v^T\off(Y_i\otimes Y_i)v)\\
&\le p^2 v^T\off(\Sigma)v + C_2'\left(\lambda_2 p^2\|\Sigma\|^2\kappa^4 + \sum_{i=1}^N \frac{R_2(Y_i)}{\lambda_2 N}+\frac{r(\Sigma)+\log(1/\delta)}{\lambda_2N}\right).
\end{split}
\end{equation*}
By Bernstein inequality the remainder terms $R_2(Y_i)$ are absorbed by the last term in the sum exactly in the same way as in the diagonal case. We optimize over $\lambda_2>0$ by setting 
\begin{equation*}
\lambda_2:= \frac{1}{p\|\Sigma\|\kappa^2}\sqrt{\frac{r(\Sigma)+\log(1/\delta)}{N}}.
\end{equation*}
Therefore, with probability $1-\delta$,
\begin{equation*}
\frac{1}{N\lambda_2}\sum_{i=1}^N \psi(\lambda_2 v^T\off(Y_i\otimes Y_i)v) \le p^2 v^T\off(\Sigma)v + C\kappa^2p\sqrt{\frac{r(\Sigma)+\log(1/\delta)}{N}}.
\end{equation*}
We repeat the arguments by changing the mean of $\nu$ to $-v$. This gives the other side of the inequality in the same way it was done for the diagonal part. We conclude that, with probability $1-\delta$,
\begin{equation*}
\|p^{-2}\widehat{\Sigma}(\lambda_2)- \off(\Sigma)\| \lesssim \frac{1}{p}\kappa^2\|\Sigma\| \sqrt{\frac{r(\Sigma)+\log(1/\delta)}{N}}.
\end{equation*}
By triangular inequality, union bound and re-scaling the multiplicative constant in $\delta$, the following holds. The estimator $\widehat{\Sigma}$ satisfies, with probability $1-\delta$, 
\begin{equation*}
\|\widehat{\Sigma}-\Sigma\| \lesssim \frac{\kappa^2}{p}\|\Sigma\|\sqrt{\frac{r(\Sigma)+\log(1/\delta)}{N}}.
\end{equation*}
\end{proof}

To end this section, we prove some technical facts, Lemma \ref{lem:almostconvex} and \ref{lemma:diagxoff}. We start with the proof of Lemma \ref{lemma:diagxoff}.
\begin{proof}
We start with the diagonal case. Observe that
\begin{equation*}
\begin{split}
&\mathbb{E}(v^T\diag(Y\otimes Y)v)^2 = \mathbb{E} \sum_{i,j=1}^d \langle Y,e_i\rangle^2\langle Y,e_j\rangle^2v_i^2v_j^2 \\
&= \sum_{i=1}^d \mathbb{E} \langle Y,e_i\rangle^4 v_i^4 + \sum_{i\neq j}^d \mathbb{E} \langle Y,e_i\rangle^2\langle Y,e_j\rangle^2v_i^2v_j^2:= (I)+(II).
\end{split}
\end{equation*}
Clearly, (I) is at most $p\kappa^4 \sum_{i=1}^d\Sigma_{ii}^2v_i^4 \le p\kappa^4\|\diag(\Sigma)\|^2$. Next, by the arithmetic-geometric inequality
\begin{equation*}
(II)\le \frac{p^2}{2}\sum_{i\neq j}\mathbb{E}(\langle X,e_i\rangle^4 v_i^2v_j^2 + \langle X,e_j\rangle^4v_i^2v_j^2) \le p^2\kappa^4\|\diag(\Sigma)\|^2.
\end{equation*}
For the off-diagonal term, we need to proceed carefully as the natural idea to decompose the off-diagonal matrix into the matrix itself minus the diagonal part leads to suboptimal dependence on $p$. We first expand it directly,
\begin{equation*}
\begin{split}
&\mathbb{E} (v^T\off(Y\otimes Y)v)^2 = \sum_{i\neq j;k\neq l}\mathbb{E}\langle Y,e_i\rangle\langle Y,e_j\rangle\langle Y,e_k\rangle\langle Y,e_l\rangle v_iv_jv_kv_l \\
&\le p^2\sum_{i\neq j;k\neq l}\mathbb{E}\langle X,e_i\rangle\langle X,e_j\rangle\langle X,e_k\rangle\langle X,e_l\rangle v_iv_jv_kv_l = p^2 \mathbb{E}(v^T\off(X\otimes X)v)^2,
\end{split}
\end{equation*}
where the term $p^2$ comes from the fact that at least two indices are distinct in each summand. Now, we split the off-diagonal term $\mathbb{E}(v^T\off(X\otimes X)v)^2$. More accurately,
\begin{equation*}
\begin{split}
&\mathbb{E}(v^T\off(X\otimes X)v)^2=\mathbb{E}(v^T(X\otimes X)v)^2\\
& + \mathbb{E}(v^T\diag(X\otimes X)v)^2 -2\mathbb{E}(v^T(X\otimes X)v)\mathbb{E}(v^T\diag(X\otimes X)v)\\
&:=(a)+(b)+(c).
\end{split}
\end{equation*}
The last term $(c)$ is negative because both matrices are positive semidefinite, so we can safely ignore it. The first term $(a)$ on is at most $\kappa^4(v^T\Sigma v)^2 \le \kappa^4\|\Sigma\|^2$ by the moment equivalence assumption. Finally, the second term $(b)$ is at most $\kappa^4 \|\Sigma\|^2$ by the same argument used above.
\end{proof}
 Next, we proceed to prove Lemma \ref{lem:almostconvex}.
\begin{proof}
Notice that $\psi(x) \le \log(1+x+x^2)$ holds trivially, and we add $x^2/6$ to make the latter function convex. It follows that  $$\psi(\mathbb{E}Z)\le \min\{\log(1+\mathbb{E}Z+\mathbb{E}Z^2) + \mathbb{E}Z^2/6,1\}.$$ Now, we apply Jensen's inequality to conclude the proof of the first part. For the second part, notice that by Taylor series expansion, if $t\in [0,a]$ then we have the following inequality,
\begin{equation*}
e^t \le 1+ \frac{t}{a}\left(\sum_{i=1}^{\infty}\frac{a^i}{i!}\right) \le 1+\frac{t}{a}(e^{a}-1),
\end{equation*}
therefore
 \begin{align*}
    &\mathbb{E}\log(1 + Z + Z^2) + a\mathbb{E}\min\{1, Z^2/6\} 
    \\
    &= \mathbb{E} \log\left(\left(1 + Z + Z^2\right)\exp(\min\{a, aZ^2/6\})\right) 
    \\
    &\le\mathbb{E} \log\left(\left(1 + Z + Z^2\right)\left(1 + \min\{1, Z^2/6\}\left(e^{a} - 1\right)\right)\right).
\end{align*}
To get the inequality in the statement, we only need to split into the cases where $|Z|^2/6$ is smaller than one and where it is greater than one. 
\end{proof}
\section{Proof of Theorem \ref{main_result}}
\label{sec:proof_main}
In the previous section, we showed in Proposition \ref{oracle_estimator}, that the proof of the main result boils down to estimate the trace of the covariance matrix, the operator norm, and the sparsifying parameter $p$. For the trace and operator norm, it is enough to estimate it with a multiplicative absolute constant. On the other hand, for the parameter $p$, we need a more accurate estimator. In fact, since we need to divide the estimator by $p$, an estimator $\widehat{p}$ that do not convergence to $p$ would insert a bias.
\begin{remark}
\label{remark_n_large}
The best possible convergence rate is at least
\begin{equation*}
   \|\widehat{\Sigma}-\Sigma\| \le \|\Sigma\|\left(\frac{1}{p}\sqrt{\frac{r(\Sigma)+\log(1/\delta)}{N}}\right).
\end{equation*}
The trivial estimator $\widehat{\Sigma}=0$ satisfies $\|\widehat{\Sigma}-\Sigma\| \le \|\Sigma\|$, so in order to have a meaningful result we need $\left(\frac{1}{p}\sqrt{\frac{r(\Sigma)+\log(1/\delta)}{N}}\right) <1$. Therefore, without making any further comments, we may assume that $$N\ge C \left(\frac{r(\Sigma) + \log(1/\delta)}{p^2}\right),$$ for some well-chosen $C>0$. 
\end{remark}
\subsection{Estimation of $p$}
The idea here is to explore the proportion of non-zeros entries in the observed data. In any standard data set, a missed value does not appear with zero; we set it to zero for convenience, as we have done throughout the entire manuscript until now. As it happens, when estimating the proportion of missing values, it could be the case that the distribution of the random vector $X$ has non-trivial mass at zero. Clearly, we can distinguish between the zero that comes from the distribution and the zero from the missing value. Equivalently, we may assume that the marginals of $X$, namely $\langle X,v\rangle$ (for every $v\in S^{d-1}$), do not have mass at zero. 

The starting point is the following. We collect $Y_1,\ldots, Y_N$, and compute $Z_1,\ldots,Z_N$, where $Z_i(j) =1$ if and only $Y_i(j) \neq 0$ and zero, otherwise. The goal is to estimate the mean of the random variable
\begin{equation*}
R(Z) := \frac{1}{d}\|Z\|_{\ell_1},
\end{equation*}
as it is equal to $\mathbb{E}R(Z)=p$.
\begin{lemma}
\label{estimation_parameter_p}
Let $Y_1,\ldots, Y_N$ be i.i.d copies of $X\odot \mathbf{p}$. There exists an estimator $\widehat{p}$ depending only on the sample and the confidence level $\delta$ satisfying that, with probability at least $1-\delta$,
\begin{equation*}
    |\widehat{p} - p| \le Cp\sqrt{\frac{\log(1/\delta)}{N}}.
\end{equation*}
As an immediate consequence, if $N\ge C\log(1/\delta)$, then (with the same probability guarantee)
\begin{equation*}
\frac{1}{2}p\le \widehat{p} \le \frac{3}{2}p.
\end{equation*}
\end{lemma}
Before we proceed to the proof, we remark that if $\widehat{p}>1$, then we round it, $\widehat{p}=1$. 
\begin{proof}
Following the notation above, we collect $R(Z_1),\ldots,R(Z_N)$ i.i.d copies of $R(Z)$. We invoke a standard sub-Gaussian mean estimator for $R(Z)$ (e.g trimmead mean estimator \cite[Theorem 1]{lugosi2021robust}) together with the fact that $\Var(R(Z))\le p^2$, to obtain that, with probability at least $1-\delta$,
\begin{equation*}
    |\widehat{R}(p) - \mathbb{E}R(Z)| \le C\sqrt{\frac{\log(1/\delta)\Var R(Z)}{N}} \le Cp\sqrt{\frac{\log(1/\delta)}{N}}.
\end{equation*}
\end{proof}

\subsection{Estimation of the Trace}
To simplify the analysis, we can safely assume that $p$ is known because we can accurately estimate it using Lemma \ref{estimation_parameter_p}. Clearly,
\begin{equation*}
    p \tr(\Sigma) = \mathbb{E} \sum_{i=1}^d \langle Y,e_i\rangle^2.
\end{equation*}
To invoke a mean estimator, we need to compute the standard deviation of the random variable in the right hand side. To this end, we have 
\begin{equation*}
    \mathbb{E} \left(\sum_{i=1}^d \langle Y,e_i\rangle^2\right)^2\le p\sum_{i=1}^d\mathbb{E}\langle X,e_i\rangle^4+p^2\sum_{i\neq j}\mathbb{E}\langle X,e_i\rangle\langle X,e_j\rangle^2 \lesssim p \kappa^4 \tr(\Sigma)^2.
\end{equation*}
The latter step follows from moment equivalence and H\"older's inequality (as we have been doing several times in this manuscript). Since $p$ is know, one may invoke Theorem 1 \cite{lugosi2021robust} to obtain an estimator $\widehat{\tr}(\Sigma)$ satisfying that, with probability $1-\delta$,
\begin{equation*}
    |\widehat{\tr}(\Sigma) - p\tr(\Sigma)| \le C \kappa^2 p \tr(\Sigma) \sqrt{\frac{\log(1/\delta)}{N}}.
\end{equation*}
If the sample size $N$ satisfies that $N\ge C\kappa^2\log(1/\delta)$ then for sufficiently large $C$, we have
\begin{equation*}
    |\widehat{\tr}(\Sigma) - p\tr(\Sigma)| \le \frac{p\tr(\Sigma)}{2},
\end{equation*}
and consequently
\begin{equation}
\label{estimation_trace}
    \frac{1}{2}\tr(\Sigma) \le p^{-1}\widehat{\tr}(\Sigma) \le \frac{3}{2}\tr(\Sigma).
\end{equation}
\subsection{Estimation of the Operator Norm}
The most delicate part of this section is the estimation of the operator norm. The main lemma is the following
\begin{lemma}
\label{estimation_operatornorm}
Let $Y_1,\ldots,Y_N$ be i.i.d copies of $X\odot p$. There exist an absolute constant $C_N$ and an estimator $\widehat{\|\Sigma\|}$ depending only on the samples and $\kappa$ satisfying that, with probability at least $1-\delta$, 
\begin{equation*}
    c_2(\kappa)\|\Sigma\|\le \widehat{\|\Sigma\|} \le c_1(\kappa) \|\Sigma\|,
\end{equation*}
provided that $N\ge C_Np^{-2}(\log(1/\delta)+r(\Sigma))$. Here $c_1,c_2>0$ are two absolute constants depending only on $\kappa$.
\end{lemma}
The key idea is to repeat the same analysis as before for each part with an additional parameter $\alpha$, and show that if certain inequalities are satisfied then $\alpha$ needs to be of same order as the operator norm.
Along the proof $C_1>0$ is an explicit constant that can be computed by just keeping track of the constants in the proofs of Section \ref{sec:oracle}.

\begin{proof}
\textbf{Diagonal Part:}
As before, we set
\[
\Theta = \mathbb{R}^d \times \mathbb{R}^d.
\]
Now, we slightly change the choice of measures. More accurately, we choose the measure $\mu$ to be a product of two zero mean multivariate Gaussians with mean zero and covariance $\beta^{-1}I_d$. For $v \in S^{d - 1}$, let $\rho_{v}$ be a product of two multivariate Gaussian distribution with mean $\alpha v$ and covariance $\beta^{-1}I_d$. The $\mathcal{KL}$-divergence becomes
\[
\mathcal{KL}(\rho_{v}, \mu) = \alpha^2 \beta.
\]
To simplify the notation, we write $\rho_{v,\alpha} = \rho_v$. Following the same lines for the proof of the diagonal part, we have with probability at least $1-3\delta$,
\begin{equation*}
\begin{split}
\frac{1}{N}\sum_{i=1}^N \psi(\alpha^2 v^T\diag(Y_i\otimes Y_i)v) &\le \alpha^2 pv^T\diag(\Sigma)v\\
&+  C_1p\|\diag(\Sigma)\|^2\kappa^4(\alpha^4+\beta^{-1}\alpha^2)\\
&+ (C_1\beta^{-2}p\kappa^4)\tr(\Sigma)^2 + \frac{2\log(1/\delta)}{N} + \frac{\alpha^2\beta}{N}.
\end{split}
\end{equation*}
Next, we choose $\beta = c_{\beta}\tr(\Sigma)$ where $c_{\beta}>0$ is an absolute constant to be chosen later. By the Remark \ref{remark_n_large}, we may define a constant $C_N>0$ for which $N\ge C_{N}p^{-2}\max\{r(\Sigma),\log(1/\delta)\}$, and then
\begin{equation*}
\begin{split}
\frac{1}{p N}\sum_{i=1}^N \psi(\alpha^2 v^T\diag(Y_i\otimes Y_i)v)& \le \alpha^2 v^T\diag(\Sigma)v +  C_1\|\Sigma\|^2\kappa^4\alpha^4\\
&+C_1\kappa^4c_{\beta}^{-1}\alpha^2\|\Sigma\|+ C_1c_{\beta}^{-2}\kappa^4 + 2C_N^{-1}\\
&+ \alpha^2\|\Sigma\|c_{\beta}C_N^{-1}.
\end{split}
\end{equation*}
\textbf{Off-Diagonal Part:} We use the same choice of the measures and proceed analogously. We obtain that, with probability at least $1-5\delta$, the following holds
\begin{equation*}
\begin{split}
&\frac{1}{p^2n}\sum_{i=1}^n \psi(\alpha^2 v^T\off(Y_i\otimes Y_i)v) \le \alpha^2 v^T\off(\Sigma)v \\
&+ C_1\alpha^4\kappa^4\|\Sigma\|^2 + C_1c_{\beta}^{-2}\kappa^4\\
&+ C_1c_{\beta}^{-1}\kappa^4\alpha^2\|\Sigma\|  + 2C_N^{-1} +C_N^{-1}c_{\beta}\alpha^2\|\Sigma\|.
\end{split}
\end{equation*}
\textbf{Everything Together:} We define the function $g(\alpha):\mathbb{R}\rightarrow \mathbb{R}$ to be equal to
\begin{equation*}
\frac{1}{Np}\sup_{v\in S^{d-1}}\sum_{i=1}^N \psi(\alpha^2 v^T\diag(Y_i\otimes Y_i)v) + \frac{1}{Np^2}\sup_{v\in S^{d-1}\cup 0}\sum_{i=1}^N \psi(\alpha^2 v^T\off(Y_i\otimes Y_i)v).
\end{equation*}
From above, we obtain that, with probability at least $1-8\delta$,
\begin{equation*}
g(\alpha) \le C_1\|\Sigma\|^2\alpha^4\kappa^4 + \|\Sigma\|\alpha^2(\kappa^2+\kappa^4C_1c_{\beta}^{-1}+c_{\beta}C_N^{-1}) + C_1\kappa^4c_{\beta}^{-2} + 4c_n^{-1}.
\end{equation*}
Notice that the constants $C_1c_{\beta}^{-1}+c_{\beta}C_N^{-1}$ and $C_1\kappa^4c_{\beta}^{-2}+4C_N^{-1}$ can be made arbitrarily small by increasing $C_N$. In particular, we choose $c_{\beta}$ and $C_N$ so that 
\begin{equation}
\label{eq:ineq_for_g}
g(\alpha) \le C_1\|\Sigma\|^2\alpha^4\kappa^4 + \|\Sigma\|\alpha^2(1+L_1) + L_2,
\end{equation}
where $L_1,L_2>0$ are two absolute constants satisfying the following conditions:
\begin{equation}
\label{eq:two_crucial_conditions_parabola}
1.1L_2<1 \quad \text{and} \quad (1-L_1)^2-8.4\kappa^4C_1L_2>0
\end{equation}
The reason for such choice will become clear in what follows. Next, without loss of generality, we assume that $\mathbb{P}(Y_i = 0)=0$. This is always possible by adding a small amount of Gaussian noise without changing the covariance too much. We construct a vector $w \in S^{d-1}$ such that $\min_{i\in [n]}\langle Y_i,w\rangle \neq 0$ by sampling the vector from an isotropic Gaussian distribution, and normalizing it to have Euclidean norm exactly one.

Notice that $g(0) =0$ and $g$ is a continuous function. Moreover, if we show that $g$ assumes values greater or equal to one then by intermediate value theorem, the function $g$ assumes any value within this range. Since $w$ is a unit vector for which $\min_{i\in [n]}\langle w,Y_i\rangle \neq 0$, and for every $i\in [N]$, $$\langle Y_i,w\rangle^2 = w^TY_i\otimes Y_iw= w^T\diag(Y_i\otimes Y_i) w + w^T\off(Y_i\otimes Y_i)w ,$$ it follows that at least one of the terms in the right-hand side is non-zero. In the case that both are non-zero, we evaluate $g$ at the point
\begin{equation}
\label{eq:min_term}
\min\left\{\min_{i\in [n]|}|w^T\diag(Y_i\otimes Y_i) w|,\min_{i\in [d]}|w^T\off(Y_i\otimes Y_i) w\}|\right\}.
\end{equation}

It is clear that $g$ is at least one at such point. Moreover, observe that the function $g$ is non-negative as we are allowed to take $v=0$ in the supremum of the off-diagonal part. In the case that one term is zero, we just remove it from \eqref{eq:min_term}. Finally, regardless the case, we choose $\widehat{\alpha}$ such that $g(\widehat{\alpha})=1.1L_2$. This is a valid choice. Indeed, recall from \eqref{eq:two_crucial_conditions_parabola} that $1.1L_2$ is strictly smaller than one, therefore existence of such $\widehat{\alpha}$ is guaranteed by the intermediate value theorem as argued before.

Next, \eqref{eq:ineq_for_g} implies that
\begin{equation*}
C_1\|\Sigma\|^2\widehat{\alpha}^4\kappa^4 + \|\Sigma\|\widehat{\alpha}^2(1+ L_1) - 0.1L_2\ge 0.
\end{equation*}
The expression above can be interpreted as a parabola in the variable $ x:=\widehat{\alpha}^2\|\Sigma\|$ that has two real roots. One root is negative, and it does not play any role. The other one is a positive absolute constant implying that there exists a  constant $c_{min}(\kappa)$ such that
\begin{equation}
\label{eq:c_min}
\widehat{\alpha}^2\|\Sigma\| \ge c_{min}(\kappa).
\end{equation}
This translates in a lower bound for $\|\Sigma\|$. We now need an upper bound for $\|\Sigma\|$ in terms of $\widehat{\alpha}$. We repeat the same argument above for the product measure $\rho_{2,v}$ between $\theta$ and $\nu$, where $\theta \sim N(\alpha v,\beta^{-1}I_d)$ and $\nu \sim N(-\alpha v,\beta^{-1}I_d)$. Therefore, if $v_1\in S^{d-1}$ is the normalized eigenvector corresponding to the maximum eigenvalue of $\Sigma$, then
\begin{equation*}
-g(\alpha) \le  -\frac{1}{np}\sum_{i=1}^n \psi(\alpha^2 v_1^T\diag(Y_i\otimes Y_i)v_1) - \frac{1}{np^2}\sum_{i=1}^n \psi(\alpha^2 v_1^T\off(Y_i\otimes Y_i)v_1).
\end{equation*}
Moreover, since $-g(\alpha)$ is non-increasing in the interval $[0,\widehat{\alpha}]$, we have
\begin{equation*}
-1.1L_2=-g(\widehat{\alpha}) \le C_1\|\Sigma\|^2\widehat{\alpha}^4\kappa^4 - \|\Sigma\|\alpha^2(1-L_1) + L_2
\end{equation*}
Setting $x= \|\Sigma\|\alpha^2$, the inequality above holds for all $\alpha \in [0,\widehat{\alpha}]$. It follows that,
\begin{equation*}
C_1\kappa^4x^2 - (1-L_1)x + 2.1L_2 \ge 0.
\end{equation*}
The discriminant of the quadratic equation is $\Delta=(1-L_1)^2-8.4C_1\kappa^4L_2$ which is (strictly) positive by \eqref{eq:two_crucial_conditions_parabola}. It follows that the inequality above is true if $x\le x_1$ or $x\ge x_2$, where $0<x_1<x_2$ are the positive roots of the corresponding quadratic equation. We claim that $x\ge x_2$ cannot happen. Otherwise, since the inequality above holds for all $\alpha \in [0,\widehat{\alpha}]$, it must hold for $\alpha^{\ast}$ such that $\|\Sigma\|\alpha^{\ast} \in (x_1,x_2)$, but this contradicts the fact that the parabola assumes negative values between $(x_1,x_2)$. Therefore, we obtain that there exists a constant $c_{max}(\kappa)>0$ such that $\widehat{\alpha}^2\|\Sigma\| \le c_{max}(\kappa)$. Putting together with \eqref{eq:c_min}, we obtain that
\begin{equation*}
   c_{min}(\kappa) \le \widehat{\alpha}^2\|\Sigma\| \le c_{max}(\kappa).
\end{equation*}
We conclude the proof  by setting $\widehat{\|\Sigma\|}:=\widehat{\alpha}^{-2}$.
\end{proof}
\subsection{Completion of the proof of Theorem \ref{main_result}}
The final construction of our estimator is the following:
\begin{tcolorbox}
\label{Box:finalEstimator}
\begin{enumerate}
    \item Split the sample $Y_1,\ldots,Y_N$ into four parts of size at least $\lfloor N/4 \rfloor$ each.
    \item Estimate the parameter $p$ with the first quarter of the sample using Lemma \ref{estimation_parameter_p}.
    \item Estimate the trace $\tr(\Sigma)$ with the second quarter using \eqref{estimation_trace} and the operator norm $\|\Sigma\|$ with the third quarter using Lemma \ref{estimation_operatornorm}.
    \item For the last quarter of the sample, use the estimator from Proposition \ref{oracle_estimator} to estimate the covariance matrix.
\end{enumerate}
\end{tcolorbox}
Before proceeding to the proof, we highlight some features about the data-splitting approach. From a theoretical perspective, it only affects the convergence rate by a constant. However, from a practical standpoint, it might be of interest to avoid wasting one quarter of the sample if there are only a few observations missing. This means that we use the complete data to estimate the covariance by setting $\widehat{p}=1$.

Moreover, it might be more appropriate to use a smaller fraction of the data to estimate the trace, as it is a one-dimensional quantity, and its convergence rate is faster than the convergence rate of the estimator itself. Unfortunately, due to the intractability of our estimator, we are unable to implement these ideas in a real dataset.

We are now in position to prove our main result, Theorem \ref{main_result}. 
\begin{proof}
As discussed in Section \ref{sec:oracle}, the proof follows easily once we estimate the parameters of the truncation level. Indeed, the truncation levels in Proposition \ref{oracle_estimator} only requires the knowledge of $\tr(\Sigma), \|\Sigma\|$ and $p$ up to an absolute constant. The error that we need to take in account is to use the estimated value of $p$ instead of the true value when we divide by $p$. This is the only reason why we have to estimate the precise value of the parameter $p$. To this end, by triangle inequality
\begin{equation*}
\left\|\frac{1}{\widehat{p}}\widehat{\Sigma}_1 -\diag(\Sigma)\right\| \le \left\|\frac{1}{\widehat{p}}\left(\widehat{\Sigma}_1 -p\diag(\Sigma)\right)\right\| + \left\|\frac{1}{\widehat{p}}\left(p\diag(\Sigma) -\widehat{p}\diag(\Sigma)\right)\right\|.
\end{equation*}
We apply Lemma \ref{estimation_parameter_p} to estimate both terms. The first term in the right hand side is, with probability at least $1-\delta$, at most 
\begin{equation*}
\frac{\sqrt{p}}{\widehat{p}} C\|\Sigma\| \sqrt{\frac{r(\Sigma)+\log(1/\delta)}{N}} \lesssim \frac{1}{\sqrt{p}}\|\Sigma \|\sqrt{\frac{r(\Sigma)+\log(1/\delta)}{N}}.
\end{equation*}
The second term also satisfies, with  probability $1-\delta$,
\begin{equation*}
\left\|\frac{1}{\widehat{p}}\left(p\diag(\Sigma) -\widehat{p}\diag(\Sigma)\right)\right\| \le \|\diag(\Sigma)\| \frac{1}{\widehat{p}}\left|p-\widehat{p}\right| \lesssim \|\Sigma\| \sqrt{\frac{\log(1/\delta)}{N}}.
\end{equation*}
The same argument holds for the off-diagonal part as clearly $\|\off(\Sigma)\|\le 2\|\Sigma\|$. We omit it for the sake of simplicity. Finally, the desired probability guarantee holds by union bound a constant number of times.
\end{proof}

\paragraph{Acknowledgments.}  
The author would like to thank Tanja Finger, Felix Kuchelmeister and Nikita Zhivotovskiy for helpful discussions.

\end{document}